\documentclass[11pt]{amsart}
\usepackage{amsmath,amsthm,amsfonts,amscd,amssymb}
\usepackage[dvips]{graphicx}

\newtheorem{thm}{Theorem}[section]
\newtheorem{lem}[thm]{Lemma}
\newtheorem{conj}[thm]{Conjecture}
\newtheorem{cor}[thm]{Corollary}
\newtheorem{prop}[thm]{Proposition}
\newtheorem{rem}[thm]{Remark}
\newtheorem{definition}[thm]{Definition}

\def\qed{$\Box$}

\def\proof{\smallbreak\noindent{\bf Proof:\ }}

\def\BB{{\bf B}}
\def\CC{{\bf C}}
\def\DD{{\bf D}}
\def\HH{{\bf H}}

\def\QQ{{\bf Q}}
\def\RR{{\bf R}}
\def\ZZ{{\bf Z}}
\def\Chat{\hat{\bf C}}
\def\Qhat{\hat{\bf Q}}
\def\Rhat{\hat{\bf R}}

\def\bch{{\partial{\cal C}}}
\def\pl{\mathop{\bf pl}}
\def\tr{\mathop{\rm Tr}}
\def\ax{\mathop{\rm Ax}}
\def\teich{\mathop{\rm  Teich}}

\def\QF{{\mathcal Q}{\mathcal F}}

\def\bch{{\partial{\mathcal C}}}
\def\pl{\mathop{ pl}}
\def\M{{\mathcal M}}
\def\DE{{\mathcal D}(\pi_1(S))}
\def\RE{{\mathcal R}(\pi_1(S))}
\def\LL{{\mathcal L}_c}
\def\DL{{\mathcal D}_c}
\def\l{\lambda}

\def\E{{\mathcal E}}
\def\F{{\mathcal F}}
\def\P{{\mathcal P}}
\def\C{{\mathcal C}}

\def\H{{\mathcal H}}
\def\Q{{\mathcal Q}}

\def\V{{\mathcal V}}

\begin{document}

\title{Linear slices of the quasifuchsian space of punctured tori}
\author{Yohei Komori}
\address{Osaka City University Advanced Mathematical Institute and
Department of Mathematics, Osaka City University, 558-8585, Osaka, Japan}
\email{komori@sci.osaka-cu.ac.jp}

\author{Yasushi Yamashita}
\address{Department of Information and Computer Sciences, Nara Women's
University, 630-8506 Nara, Japan}
\email{yamasita@ics.nara-wu.ac.jp}
\date{}

\begin{abstract}
After fixing a marking $(V, W)$ of a quasifuchsian punctured torus
group $G$, the complex length $\l_V$ and the complex twist
$\tau_{V,W}$ parameters define a holomorphic embedding of the
quasifuchsian space $\QF$ of punctured tori into $\CC^2$. It is called
the complex Fenchel-Nielsen coordinates of $\QF$.  For $c \in \CC$,
let $\Q_{\gamma, c}$ be the affine subspace of $\CC^2$ defined by the
linear equation $\l_V=c$.  Then we can consider the {\em linear slice}
$\LL$ of $\QF$ by $\QF \cap \Q_{\gamma, c}$ which is a holomorphic
slice of $\QF$.  For any positive real value $c$, $\LL$ always
contains the so called {\em Bers-Maskit slice} ${\mathcal BM}_{\gamma,
  c}$ defined in~\cite{KS}.  In this paper we show that if $c$ is
sufficiently small, then $\LL$ coincides with ${\mathcal BM}_{\gamma,
  c}$ whereas $\LL$ has other components besides ${\mathcal
  BM}_{\gamma, c}$ when $c$ is sufficiently large.  We also observe
the scaling property of $\LL$.
\end{abstract}

\maketitle

\section{Introduction}
The quasifuchsian space $\QF$ of once punctured tori can be embedded in
$\CC^2 = \{ (\lambda, \tau) \}$ by the complex Fenchel-Nielsen coordinates
(c.f. \cite{KS,  Christos, PP, SerTan}).
By varying the complex twist $\tau$ and keeping the complex length $\lambda$
being fixed as
a positive real value $c$,
we can define the {\em linear slice} $\LL \subset \CC$ of $\QF$.
In this paper we investigate the global
properties of $\LL$ realized in the complex plane.  To state our results, recall
that $\LL$ has a component containing the open interval $(2, +\infty)$
which was studied in~\cite{KS} and also  in~\cite{KoP, PP}.  In this paper we call this
component the {\em standard component} and the others {\em non-standard}.
We will show
\setcounter{section}{5}\setcounter{thm}{0}
\begin{thm}
There exists some positive constant $c_0$ such that for any $c$
satisfying $0 < c < c_0$, $\LL$ coincides with the standard component.
\end{thm}
\setcounter{section}{6}\setcounter{thm}{0}
\begin{thm}
There exists some positive constant $c_1$ such that for any $c$
satisfying $c >c_1$, $\LL$ contains non-standard components.
\end{thm}
In section 7, we also consider the scaling property of $\LL$.
\setcounter{section}{7}\setcounter{thm}{2}
\begin{cor}
Linear slice has an asymptotic scaling constant.
\end{cor}

See Figure 1 for theorem 5.1 and Figure 2 and 3 for theorem 6.1 and
corollary 7.3.  The parameters used in the figures are explained in
4.1.

\setcounter{section}{1}\setcounter{thm}{0}

Let us describe some historical background of our subject.
A marked quasifuchsian punctured torus group $G$ is a free marked two generator
discrete subgroup of $PSL_2({\bf C})$ such that the commutator of the
generators
is parabolic, and the regular set $\Omega$ consists of two non-empty simply
connected
invariant components $\Omega^{\pm}$.
Quasifuchsian space $\QF$ is the space of marked quasifuchsian punctured
torus groups
modulo conjugation in $PSL_2({\bf C})$.
The convex core ${\mathcal C}/G$ has two boundary components $\partial
{\mathcal C}^{\pm}/G$
each of which is a once-punctured torus and admits an intrinsic hyperbolic
structure
making it a pleated surface.

In their seminal paper~\cite{KS} L.~Keen and C.~Series
defined the
{\em Bers-Maskit slice} ${\mathcal BM}_{\mu, c}$ 
for a fixed measured lamination $\mu$ and $c>0$, 
as the subset of $\QF$ on
which
the bending lamination of
$\partial {\mathcal C}^-/G$ and $\mu$ belong to the same projective class
and the length of $\mu$ in $\partial {\mathcal C}^-/G$ is equal to $c$.
By using their theory of pleating coordinates, they showed that ${\mathcal
BM}_{\mu, c}$ is simply
connected.
J.~Parker and J.~Parkkonen also studied these slices for the case that
$\mu$ is a
rational lamination
(they call them the $\lambda$-{\em slices}), and considered a generalization
of I.~Kra's
plumbing construction and degeneration of ${\mathcal BM}_{\mu, c}$ 
to the Maskit slice ${\mathcal M}$
(c.f.~\cite{PP}).
The first author and J.~Parkkonen  further studied ${\mathcal BM}_{\mu, c}$;
they showed that the boundary of ${\mathcal BM}_{\mu, c}$ is a Jordan curve which is cusped at a countable dense set of points
(c.f.~\cite{KoP}).
In this papaer
we would like to study the
outside of ${\mathcal BM}_{\mu, c}$  in $\LL$ and  its scaling property.

This paper is organized as follows.  In section 2 
we will review the basic notions of the quasifuchsian space $\QF$ of once punctured tori 
and its pleating varieties following~\cite{KS}.
The complex Fenchel-Nielsen coordinates of $\QF$ will be introduced in section 3, and 
we will define the main subject of this paper, the {\em linear
slice} $\LL$ of $\QF$  in section 4.
In sections 5 and 6 we will study connected components of $\LL$ and prove our main theorems.
And in the last section we will observe
the asymptotic self similarity
of $\LL$.


We are grateful to Caroline Series for showing us the preprint~\cite{Otal}
of Otal,
and Raquel Diaz for explaining us her idea in section 6.
It was fruitful for us to discuss with them in Nara in January 2000;
in practice this work was almost done during their stay in Japan.
We also wish to thank Hideki Miyachi for enjoyable conversations with him on the
topic in section 7, and Kentaro Ito and Sara Maloni for telling  their interests in our paper recently.

The first author was partially supported by Grant-in-Aid for Scientific Research(C) (19540194), Ministry of Education, Science and Culture of Japan.

\section{The quasifuchsian space $\QF$ and rational pleating varieties}
\subsection{Punctured torus groups and their pleating data}
\label{sec:normalization}
\subsubsection{Marking}
Let $S$ be an oriented once-punctured torus.
Any pair of simple closed loops on $S$
that intersect exactly once are free generators of
$\pi_1(S)$.
 Let $(\alpha,\beta)$ be
such an ordered pair of free generators, chosen so that their
commutator
$\alpha \beta \alpha ^{-1} \beta ^{-1}$ represents a  positively
oriented loop around the puncture.
  The ordered pair $(\alpha,\beta)$ is called
a {\em marking}.

\subsubsection{$\QF$ and $\F$}
A {\em punctured torus group} is a discrete subgroup
$G \subset  PSL(2,\CC)$ that  is
the image of a faithful
representation $\rho$ of  $\pi_1(S)$
such that the image of the loop around the puncture is parabolic.
If  $(\alpha,\beta)$ is a   marking of $S$,
and if  $
A = \rho(\alpha),B =\rho(\beta)$, then
the commutator  $K = ABA^{-1}B^{-1}$
is parabolic and  the ordered pair
$(A ,B) = (\rho(\alpha),\rho(\beta))$ is called a
{\em marking} of $G$.

The group $G$ is {\em quasifuchsian} if the regular set $\Omega(G)$
consists
of two non-empty simply connected invariant components
$\Omega^{\pm}(G)$.
The limit set $\Lambda(G)$ is topologically a circle.
{\em Quasifuchsian space} $\QF$ is the space of marked
quasifuchsian
punctured torus groups modulo conjugation in $PSL(2,\CC)$; it has
a holomorphic structure induced from the
natural holomorphic structure of $SL(2,\CC)$.

Let $\RE$ be the set of $PSL(2,\CC)$-conjugacy classes of representations
$\rho$ of  $\pi_1(S)$
such that the image of the loop around the puncture is parabolic.
Considering the compact open topology on $\RE$,
Minsky showed that the closure of $\QF$ in $\RE$ is equal to $\DE$,
the set of punctured torus groups modulo conjugation in $PSL(2,\CC)$
(c.f.~\cite{Min}).

{\em Fuchsian
space} $\F$ is the subset of $\QF$ such that the components $\Omega^{\pm}$
are round
disks.  It is canonically isomorphic to the Teichm\"uller space of
marked
conformal structures on
$S$.

The quotients
$\Omega^{\pm}(G)/G$ are punctured tori with conformal structures, and
hence
also with orientations inherited from $\Chat$; we assume that the orientations of
$\Omega^+(G)/G$
and $S$ agree whereas those of $\Omega^-(G)/G$ and $S$ are
opposite.

A point $q \in \QF$ represents an equivalence class of marked
groups
in $ PSL(2,\CC)$.
We choose once and for all a  triple of distinct
points in $\Chat$ and let $G= G(q)$ denote the
representative normalized so that the repelling and attracting
fixed points of $A$ and the fixed point of $K$
are equal to the fixed triple points in $\Chat$.
If it is clear from the
context, for readability, we suppress the dependence on $q$.

\subsection{Simple closed curves}
\label{sec:Simple closed curves}

\subsubsection{Enumeration}
\label{sec:Enumeration}
Denote by $C(S)$, the set of free unoriented homotopy
classes of
simple closed non-boundary
parallel curves on
$S$.
As is well known, this set can be naturally identified with
$\Qhat =\QQ \cup \infty$.
One way to see this is as follows;
Let ${\mathbf L}$ denote the integer lattice $m+in, m,n \in \ZZ
\subset
\CC$.
Topologically $S$ is the
quotient of the punctured plane $\CC_{i}= \CC-{\mathbf L}$
by the natural action of
$G_{i}= \langle \hat A, \hat B_{i} \rangle \equiv \ZZ^2$
by the horizontal and vertical translations.
A straight line  of rational slope
in $\CC-{\mathbf L}$ projects onto a
simple closed curve on the marked punctured torus
$S_{i} = \CC_{i}/G_{i}$, and the projection of all lines
of the same rational slope with the same orientation are
homotopic.
We denote the unoriented homotopy class obtained by projecting
the
line of
slope
$-q/p$ by $[L(p/q)]$.
Relative to our choice of
marking,
$[L(p/q)]$
is in the homology class of
$\alpha^{-p}\beta^q$ or $\alpha^{p}\beta^{-q}$ on $S_{i}$,
where $\alpha, \beta$ are projections
of the horizontal and vertical lines corresponding to $\hat A, \hat
B$
respectively.
Setting $1/0 = \infty$, we obtain that the map
${\Qhat} \rightarrow C(S)$ defined by
$p/q \mapsto [L(p/q)]$ which is well-defined and bijective.
The reason for the choice of convention
that $[L(p/q)]$ corresponds to
$\alpha^{-p}\beta^q$,  is the following; if we identify
the Teichm\"{u}ller space $\teich(S)$ of once punctured
tori
with the upper half plane ${\bf H}$,
then one can easily compute that the boundary point $p/q \in
{\bf
\Rhat}$
is the point where the extremal length of curves in the class
$[L(p/q)]$ has shrunk to zero.

\subsubsection{Special word $W_{p/q}$}
\label{sec:special word}
Suppose that $\rho: \pi_1(S) \to G \subset PSL_2(\CC)$
is a quasifuchsian punctured
torus group, marked as usual by generators $A = \rho(\alpha), B
=
\rho(\beta)$.
We denote the unique geodesic in the homotopy class of
$\rho([L(p/q)])$ in $\HH^3/ G$
by $\gamma_{p/q}$. In particular, for $q \in \QF$, 
$\gamma_{p/q}(q)$ represents the corresponding geodesic in $\HH^3/
G(q)$.

For each $p/q \in \Qhat$, we can find an
explicit word  $W_{p/q}$ in the marked generators $\langle
\alpha,
\beta
\rangle$
of $\pi_1(S)$
representing  $[L(p/q)]$ as follows.
The  words are generated
from the initial data
$$
W_{0/1}=\beta, \; W_{1/0}=\alpha^{-1}
$$

by the formula

$$
W_{(p+r)/(q+s)}=W_{r/s}W_{p/q}
,$$
whenever   $p/q<r/s$ and  $ps-qr=- 1$.
We denote by
$W_{p/q}(q)$ the corresponding special word in $G(q)$.

\subsection{Rational pleating varieties}

\subsubsection{The pleating loci}
\label{sec:pleating locus}

We are now ready to discuss the convex hull boundary and the
pleating
locus.
Let $q \in \QF$ and let   $G=G(q)$
be the corresponding  marked quasifuchsian group with the regular
set and the limit set
$\Omega(G),\Lambda(G)$ respectively.
The $3$-manifold  $ \HH^3/G$ is homeomorphic to
$S \times (0,1)$. The
surfaces
$\Omega(G)/G$ at infinity
form the boundary $S \times \{0,1\}$.
Let $\bch (G)$
be the boundary  of the
hyperbolic convex hull of $\Lambda(G)$ in $\HH^3$;
it is clearly invariant under the action of $G$.
The nearest point retraction
$  \Omega(G) \rightarrow \bch (G)$,
defined by mapping $x \in \Omega(G)$ to
the unique point of contact with $\bch (G)$ of the largest
horoball in
$\HH^3$
centered at $x$  with interior disjoint from $\bch (G)$, can
easily be
modified to
a $G$-equivariant
homeomorphism.
We denote
two connected components of $\bch (G)$
corresponding to $ \Omega^{\pm} (G)$ by $\bch^{\pm} (G)$
respectively.
Thus each component $ \bch^{\pm} (G)/ G$ is topologically
a  punctured torus. (In the special case in which $G$ is
Fuchsian,
 $\bch (G)$ is a flat plane  whose two sides serve as a
substitute for
the two components $ \bch^{\pm} (G)$.)

$ \bch^{\pm} (G)/G$ are pleated surfaces in
$\HH^3/ G$.
More precisely, there are complete hyperbolic
surfaces $S^{\pm}$, each homeomorphic to $S$,
and  maps $f^{\pm} : S^{\pm} \rightarrow \HH^3/G$,
such that
every point in $ S^{\pm}$ is in the interior of some
geodesic arc
which is
mapped by $f^{\pm}$ to a geodesic arc
in $\HH^3/G$, and such that
$f^{\pm}$ induce isomorphisms $\pi_1(S) \to G$.
Further, $f^{\pm}$ are isometries onto their images with the
path
metric induced from $\HH^3$ (c.f.~\cite{Thur}).
The {\em bending} or {\em pleating locus} of $
\bch^{\pm} (G)/G$
consists of those points of $S^{\pm}$ contained in the
interior
of one and only one geodesic arc which is mapped by $f^{\pm}$ to
a
geodesic arc in $\HH^3/G$.
For $G$ non-Fuchsian, the pleating loci are
geodesic laminations, meaning that they are unions of pairwise
disjoint
simple
closed geodesics on $ S^{\pm}$.  We denote these laminations
by
$|\pl^{\pm}(q)|$, and usually identify such a lamination with
its
image
under $f^{\pm}$ in
$\HH^3/G$.  A geodesic lamination is called {\em
rational}
if it consists entirely of closed leaves.
We concentrate on the special case in which
at least one of the pleating loci is rational in this sense.
Since
the maximum number of pairwise
disjoint simple closed curves on a punctured torus is one, such
a
lamination consists of a single
simple closed geodesic and is
therefore of the form $\gamma_{p/q}(q)$ for some
$p/q \in \Qhat$.

\subsubsection{Rational pleating varieties and hyperbolic loci}
\label{sec:plvars}

Given $p/q \in \Qhat$, we set
$$\textstyle \P_{p/q}^{\pm}= \{ q \in \QF : |pl^{\pm}(q)|= \gamma_{p/q}(q) \}
\mbox{ and }
\P_{p/q}= \P_{p/q}^+ \cup \P_{p/q}^-.$$
We call these sets the $p/q$-{\em pleating varieties}.

For any $p/q \in \Qhat$,
consider the trace
$\tr W_{p/q}$ of the special word $W_{p/q}$ associated to
$p/q$ defined in~\ref{sec:special word}.
For $ q \in \RE$,
we may consider
the function $T_{p/q}(q) = \tr W_{p/q}(q)$ as a
rational function on  $\RE$.
We define the
{\em   hyperbolic locus} of $T_{p/q}$
to be the set
$$
\H_{p/q}=\{ q \in \RE : T_{p/q}(q)\in \RR, \;
|T_{p/q}(q)|>2  \}.
$$

Then the next result is fundamental (c.f. Proposition 22 in~\cite{KS}).
\begin{prop}
$\P_{p/q} \subset \H_{p/q}$.
\qed
\end{prop}

\section{The complex Fenchel-Nielsen coordinates of $\QF$}

\subsection{The complex length of a loxodromic element }
\label{sec:cxlenlox}

The {\em complex translation length}
$\lambda_M \in \CC/2 \pi i \ZZ$ of $M \in PSL(2,\CC)$
 is given
by the equation
\begin{equation}
\label{eqn:trmult} {\pm} \tr M = 2 \cosh {\lambda_M/2}
,\end{equation}
where $\tr M$ is the trace of $M$ and we
choose the sign so that $\Re \lambda_M \ge 0$.

The complex length is invariant under conjugation by M\"obius
transformations
and has the following geometric interpretation, provided $M$ is not
parabolic; let $x$ be a point on the axis $\ax M$ of $M$ and let $\bar v$ be a vector normal
to $\ax
M$ at $x$.  Then $\Re \lambda_M$ is the hyperbolic distance between
$x$ and
$M(x)$ and $ \Im \lambda_M$ is the angle ${\rm mod}\; 2 \pi$ between
$M(\bar v)$
and the parallel transport of $\bar v$ to $M(x)$, measured facing
the
attracting fixed point $M^+$ of $M$.  In particular, if $M$ is
loxodromic
then $ \Re \lambda_M > 0$ and if $M$ is purely hyperbolic then in
addition
$ \Im \lambda_M \in 2\pi\ZZ$.

For  $q \in \QF$ and $\gamma \in C(S)$,
we denote the
element  in the group $G(q)$ representing $\gamma$ by
$W(q)$.  Because the trace is a conjugation invariant,
the complex translation length $\lambda_W(q)$ depends only on $q$
and is
independent of the normalization
of $G(q)$.
We want to define the complex length
$\lambda_{\gamma}(q)=\lambda_W(q)$ as a holomorphic
function on $\QF$ with values in
$\CC$, not $\CC/2 \pi i \ZZ$.
To do this, we choose the branch that is real valued on $\F$.
Since $\lambda_{\gamma}
\neq 0$ on $\QF$ this choice uniquely determines a  holomorphic
function
$\lambda_{\gamma}:\QF
\rightarrow \CC$.
From now on, the term
``complex length'' will always refer to this  branch.

We remark that $\Re \lambda_{\gamma}(q)$ is the hyperbolic length of $\gamma$
in $\HH^3/G(q)$.

\subsection{The complex Fenchel-Nielsen coordinates}
\label{sec:cfncoords}

The complex Fenchel-Nielsen
parameters were introduced in \cite{Christos,SerTan}
as a generalization to $\QF$ of the classical
Fenchel-Nielsen coordinates for Fuchsian groups.  Here we briefly
summarize
the main points as applied to the case of a punctured torus $S$.

Let $G=\langle  A,B \rangle$ be a marked quasifuchsian punctured
torus
group constructed from a pair of marked generators $\alpha, \beta$
of
$\pi_1(S)$ as described in~\ref{sec:normalization}.  
The complex Fenchel-Nielsen coordinates $(\lambda_A, \tau_{A,B})$ for $G= \langle  A,B
\rangle$
are obtained as follows; the parameter $\lambda_A \in \CC /2\pi i
\ZZ$ is
the complex translation length of the generator $A=\rho(\alpha)$,
or
equivalently the complex length $\lambda_{\alpha}$. The twist
parameter
$\tau_{A,B} \in \CC /2\pi i \ZZ$ measures the complex shear when
the axis
$\ax B^{-1}AB$ is identified with the axis $\ax A$ by $B$.  More
precisely,
if the common perpendicular $\delta$ to $\ax B^{-1}AB$ and $\ax A$
meets
these axes in points $Y,X$ respectively, then $\Re \tau_{A,B}$ is
the {\em
signed} distance from $X$ to $B(Y)$ and $\Im \tau_{A,B}$ is the
angle
between $\delta$ and the parallel translate of $B(\delta)$ along
$\ax A$ to
$X$, measured facing towards the attracting fixed point of $A$.
The conventions for measuring the signed distance and the angle are
explained
in more detail in \cite{KSbend}.

As shown in ~\cite{PP, Christos, KSbend}, given the parameters
$ \lambda_A,\tau_{A,B}$, and fixed a normalization, one can
explicitly
write down the matrix generators for a marked
two generator group $ G(\lambda_A,\tau_{A,B}) \subset PSL(2,\CC)$
in which the commutator $[A,B]$ is
parabolic as follows:

$$
A  =  \left(
\begin{array}{cc}
\cosh(\frac{\lambda}{2}) & \cosh(\frac{\lambda}{2}) +1 \\
\cosh(\frac{\lambda}{2}) -1 & \cosh(\frac{\lambda}{2})
\end{array}
\right),
$$
$$
B  =  \left(
\begin{array}{cc}
\cosh(\frac{\tau}{2}) \coth(\frac{\lambda}{4}) & -\sinh(\frac{\tau}{2})  \\
-\sinh(\frac{\tau}{2})  & \cosh(\frac{\tau}{2}) \tanh(\frac{\lambda}{4})
\end{array}
\right).
$$

 This group may or may not be discrete.  The matrix
coefficients
of $G$ depend holomorphically on the parameters.  The
construction thus defines a holomorphic embedding of $ \QF$ into a
subset
of $\CC /2\pi i \ZZ \times \CC /2\pi i \ZZ$, in which Fuchsian
space $\F$
is identified with the image of $\RR^2$.

We want to lift this to an embedding into $\CC^2$.
In~\ref{sec:cxlenlox} we  discussed  how to lift the length
function
$\lambda_A$ on $\QF$ to a holomorphic function  on  $\CC$.
We can similarly lift the twist parameter $\tau_{A,B}$ by
specifying that
it  will be real valued on $\F$.

On  $\F$, the  real valued parameters $ \lambda_A,\tau_{A,B}$
reduce to the classical Fenchel-Nielsen
parameters $ l_A,t_{A,B}$ defined by the above construction with
 $\lambda_A$  the
hyperbolic translation length $ l_A$ of $A$ and
$\tau_{A,B}$  the twist parameter $t_{A,B}$.

\subsection{Rational quakebends and pleated surfaces}
\label{sec:genqbs}

Clearly, the complex Fenchel-Nielsen coordinates can be made
relative to any
marking $V,W$  of $G$.
As described in detail in section~5 of \cite{KSbend}, for fixed
$\lambda \in \RR^+$ and $\tau \in \CC$,
the complex Fenchel-Nielsen coordinates relative to
$V,W$
determines a pleated surface $\psi: \DD \to \HH^3$.
We review this process.

Write $\V$ for the set of all lifts of the simple closed curve $\gamma$
corresponding to $V$ to $\DD$.
Since $\gamma$ is simple, $\V$ consists of a set of pairwise disjoint
geodesics in $\DD$, namely the axis of $V$ and all of its conjugates.
These axes in $\V$ partition $\DD$ into pieces $P_i$.
The map $\psi$ is defined in such a way that $\psi$ is an isometry on each
axis in $\V$
and on each closed piece $P_i$.
Let $x,y \in \DD - \V$ and let $\beta$ be an oriented geodesic from $x$ to $y$.
Let $P_0, P_1, ... , P_k$ be the pieces cut in order by $\beta$, that meet
along axes
$\alpha_1, \alpha_2, ... , \alpha_k \in \V$.
Orient $\alpha_i$ so that, in $\DD$, $P_{i-1}$ lies to the left of
$\alpha_i$ and $P_i$
to the right. Let $X_i=\beta \cap \alpha_i$ and
let $\bar{v}_i, \bar{w}_i$ be tangent vectors to $\psi(P_{i-1} \cap \beta)$ and
$\psi(P_i \cap \beta)$ at $\psi(X_i)$,
oriented in the direction inherited from $\beta$, so that $\bar{v}_i$ points
out of $\psi(P_{i-1})$
and $\bar{w}_i$ points into $\psi(P_i)$.
Let $\bar{v}'_i, \bar{w}'_i$ be the projections of $\bar{v}_i, \bar{w}_i$
onto the directions
orthogonal to the image of the bending axis at $\psi(X_i)$.
Then $\Im \tau$ is the angle from $\bar{v}_i$ to $\bar{w}_i$ measured facing
along $\psi(\alpha_i)$.
We embed $\DD$ in the hyperbolic ball model $\BB^3$ of $\HH^3$
as the equatorial plane such that the origins in $\DD$ and in $\BB^3$ coincide.
We arrange that the axes of $V$ and $WVW^{-1}$ in $G(\lambda, \Re \tau)$ lie
in the
boundary of a piece $P_0$ contained in $\DD$.
We then choose $\psi$ to be the identity on $P_0$.
We set $ \DD_{\gamma}(\lambda, \tau)= \psi(\DD) $
for the image of the pleated surface in $\BB^3$.
Then $\psi$ induces the group isomorphism
$\psi_*:G(\lambda, \Re \tau) \rightarrow G(\lambda, \tau)$
satisfying that
$\psi(g(z))=\psi_*(g)(\psi(z))$ for $g \in G(\lambda, \Re \tau)$ and $z \in
\DD$.

The next proposition explains the relation between $\psi$ and the bending locus
of $\bch^-(G(q))$
for $q \in \QF$.

\begin{prop}
\label{prop: homeo}
For $q \in \QF$, let $(\lambda, \tau)$ be the complex Fenchel-Nielsen
coordinates
relative to marked  generators $(V,W)$ of $G(q)$,
and let $\gamma$ be the simple closed curve
corresponding to $V$.
Assume that $V$ is purely hyperbolic and let $\psi: \DD \to \HH^3$ be
the pleated surface defined above. Then  $\psi$ is a
homeomorphism if and only if $|\pl^-(q)|= \gamma$.
\end{prop}

\proof
First suppose that $\psi$ is a homeomorphism.
Then the boundary of
$ \DD_{\gamma}(\lambda, \tau)$ is $\Lambda(G(q))$.
$ \DD_{\gamma}(\lambda, \tau)$ divides $\HH^3$ into two domains;
one of which is convex, hence contains $\C(G(q))$.
Moreover $ \DD_{\gamma}(\lambda, \tau)$ contains the axis of $V$ and all of
its conjugates in $G(q)$,
and the complement of them consist of totally geodesic pieces.
Therefore it is one of the component of $\bch(G(q))$ and from the bending
construction
in the above argument, it should be equal to $\bch^-(G(q))$
(c.f. section 7.1 in~\cite{KSbend}).

Next suppose that $|\pl^-(q)|= \gamma$.
Then $\bch^-(G(q))$ consists of the axis of $V$ and all of its conjugates in
$G(q)$,
and totally geodesic pieces. The stabilizer subgroup of each totally
geodesic piece is conjugate
to the Fuchsian subgroup $\langle V, WVW^{-1} \rangle$.
Therefore we can construct the pleated surface satisfying that $
\DD_{\gamma}(\lambda, \tau)= \bch^-(G(q))$,
which implies that $\psi$ is a homeomorphism.
\qed

\subsection{Rational quakebend planes}
\label{sec:qbplanes}
Let $(\lambda_V,\tau_{V,W})
\subset
\CC^2$ be the complex Fenchel-Nielsen coordinates relative to
marked  generators $(V,W)$ of $G$,
and let $\gamma$ be the simple closed curve
corresponding to $V$.
Assume that $V$ is purely hyperbolic and let $c$ be the hyperbolic length of
$\gamma$
in $\HH^3/G$.

We denote the slice $ \{(c,\tau) \in \CC^2 \; | \tau \in \CC \}$
by $\Q_{\gamma, c}$ and call it the {\em rational quakebend plane}.

Clearly,
$\Q_{\gamma, c}$  meets $\F$ along the earthquake path
(c.f.~\cite{KS}).
The quakebend parameter
$\tau$ is a holomorphic coordinate on $\Q_{\gamma, c}$.

On $\Q_{\gamma, c}$, the quakebend parameter
$\tau$ and $\tr W$ are related by
$$
\tr W = 2 \coth(\frac{c}{2}) \cosh(\frac{\tau}{2}).
$$
On $\Q_{\gamma, c}$, $\tr W$ is a holomorphic function of $\tau$,
branched at $\tau = 2\pi i n \;
(n \in \ZZ)$.
(see figure 5.1 in~\cite{PP}).
When $\tr V$ is real, $\QF \cap \Q_{\gamma, c}$ is contained in the strip
$$
\{ \tau \in \CC \; | \;  -\pi i < \Im \tau < \pi i \}
$$
from the argument in~\ref{sec:genqbs}.
$\tr W$ takes the right half strip
$$
\{ \tau \in \CC \; | \; \Re \tau >0, \;  -\pi i < \Im \tau < \pi i \}
$$
conformally onto
the right half plane $\CC^+$ minus the interval $( 0 ,2 \coth(\frac{c}{2})]$
where the interval $\{ \tau \in \CC \; | \; \Re \tau = 0, \; -1< \Im \tau <
1 \}$
in the imaginary axis
is folded at the origin by $\tr W$ and its image is $( 0 ,2
\coth(\frac{c}{2})]$.
We remark that $\QF \cap \Q_{\gamma, c}$ is also periodic under the action
of the Dehn twist $(A, B) \mapsto (A, A^n B)$,
and symmetric under the holomorphic involution $\tau \mapsto -\tau$.

\section{The linear slice $\LL$}
\subsection{Definition}
\label{sec:linear slice}
For $q \in \RE$, a marked group $G(q)=\langle A,B \rangle$ modulo
conjugation in $PSL(2,\CC)$
is uniquely determined by $\tr A, \tr B$ and $\tr AB$.
In fact, ignoring marking, $G(q)$ modulo conjugation in $PSL(2,\CC)$
is determined only by $\tr A$ and $\tr B$
(more precisely, the pair $(\tr A, \tr B)$ determines a marked group
$\langle A,B \rangle$ or $\langle A,B^{-1} \rangle$
modulo conjugation in $PSL(2,\CC)$).
As an application of the Jorgensen's theory on the combinatorial structure
of the Ford domain of a punctured torus group,
there is an algorithm roughly answering whether $G(q)$ is a
geometrically finite discrete group or not
from the data $(\tr A, \tr B)$ (c.f.~\cite{ASWY}).
Especially fixing $\tr A=c$, then we can use this algorithm to draw the
picture of
$$
\DL= \{ \tr B \in \CC^+ \; | \; G(q)=\langle A,B \rangle \; \text{is a
geometrically finite discrete group} \}.
$$
We call this set the {\em discrete locus}.
Let $\Q_{\gamma, c}$ be the rational quakebend plane
in the complex Fenchel-Nielsen coordinates relative to
the corresponding marked  generators $(A,B)$ of $G$ where we assume that
$c = \lambda_{A}(q)$ is real.
The {\em linear slice} $\LL$ in the right half plane $\CC^+$, which is the image of
$\QF \cap \Q_{\gamma, c}$ under $\tr B$.
Because $\QF$ is open in $\CC^2$ in complex Fenchel-Nielsen coordinates
and $\tr B$ is a open map on $\Q_{\gamma, c}$, $\LL$ is open in $\CC^+$.
Then from the definition of $\DL$ and $\LL$, $\LL$ is a subset of the
interior of $\DL$.

\begin{prop}
The interior of $\DL$ is equal to $\LL$.
\end{prop}

\proof
It is enough to show that any point $q_0$ of $\DL$,  not contained in $\LL$
is a boundary point of $\DL$.
First suppose that $G(q_0)$ is not a free group. Then some word, say $g(q_0)$
is trivial in $G(q_0)$.
Then applying the Jorgensen's inequality for the subgroup $H(q_0)$ generated by
$g$ and $K=[A, B]$,
we can see that if we take a small neighborhood $U$ of $q_0$,
for any point $q$ of $U$ except $q_0$,
$H(q)$ is not discrete which means that $G(q)$ is also
indiscrete. Therefore $q_0$ is an isolated point of $\DL$.
Next suppose that $G(q_0)$ is free, hence geometrically finite
non-quasifuchsian punctured torus group.
Then it must be a cusp (c.f.~\cite{Min}). Hence there is some word, say
$g(q)$ which is parabolic in $G(q_0)$.
Since $\tr g$ is a holomorphic function of $\tr B$, it is a open map, hence
there is
a path in the $\tr B$-plane starting from $q_0$
such that on this path $g$ is elliptic. Therefore this path is outside of
$\DL$.
This implies that $q_0$ is a boundary point of $\DL$.
\qed

\vspace{3mm}

From this result, we can see $\LL$ as the interior $\DL$
and study them experimentally.
Figures at the end of this paper drown by the second author show
computer-generated linear slices, revealing some global properties.
The black region corresponds to the discrete locus.
In the first picture, $\tr A$ is fixed at $2$ and $\tr B$ ranges in
the square of width $4$ centered at $\tr B = 2$ so that we see the
familiar picture of the Maskit slice.
By setting $\tr A = 2.5$, we get the second picture.
In Figure 2 and 3,
the value of $\tr A$ is fixed at
$8$ and $100$ respectively while changing the ranges of $\tr B$.
The width of the squares are 16, 32, 128 and 128, 2560, 12800
respectively.
We can clearly see the ``rough self similarity'' of the pictures
between figures 3 and 5 and between figures 6 and 8,
with which we will discuss in section 7.

\subsection{Connected components of $\LL$}
\label{sec:components}

\begin{prop}
For any $c>0$, $\LL$ has a component containing an open interval $(2, +\infty)$.
\end{prop}

\proof
There exists a component in $\QF \cap \Q_{\gamma, c}$ containing $\F \cap \Q_{\gamma, c}$ the real line
which is
periodic under the action of the Dehn twist
$B \mapsto A^nB$ and symmetric under $\tau \mapsto -\tau$
(c.f.~\cite{KS, PP}. see also~\cite{Mc}).
Then its image under $\tr B$ is the required component.
\qed

This component is called the {\em BM-slice} in~\cite{KS}
and also called the {\em $\lambda$-slice} in~\cite{PP}.
In this paper we call this component of $\LL$
the {\em standard component}, and call the other components
the {\em non-standard components} if they exist.
Because of the existence of the standard component which contains the
critical value of $\tr B$,
if there is a non-standard component, it is a conformal image of a component
of $\QF \cap \Q_{\gamma, c}$
under the map $\tr B$. Therefore
we can consider that $\LL$ describes the picture of $\QF \cap \Q_{\gamma, c}$.
Next result shows that topologically every component is a disk.

\begin{prop}
Each component of $\LL$ is simply connected.
\end{prop}
\proof
This is a consequence of a result of McMullen~\cite{Mc} that
$\QF$ is disk convex in $\RE$,
that is,
for any continuous map $f: \overline{\Delta} \rightarrow \RE$
whose restriction to the unit disk $\Delta$ is holomorphic,
$f(\partial \overline{\Delta}) \subset \QF$ implies
$f(\overline{\Delta}) \subset \QF$.
\qed

\begin{rem}{\rm (The Maskit slice)}\\
If we consider the limiting case where $c=0$,
we can no longer consider the complex Fenchel-Nielsen coordinates.
But by using $\tr B$ we can realize the part of the boundary of $\QF$ defined by
the
condition that $A$ is parabolic.
Then the standard component defined above corresponds to the so-called
Maskit slice $\M$.
\end{rem}

\section{Non-existence of non-standard components}
\subsection{Otal's result}
\label{sec: Otal}
\begin{thm}
\label{thm: non-existence}
There is some positive constant $c_0$ such that for any $c$ satisfying $0< c
< c_0$,
$\LL$ coincides with the standard component.
\end{thm}

This is an immediate consequence of the following result due to
J.~P.~Otal~\cite{Otal}.

\begin{thm}{\rm (c.f. corollaire 9.1 in~\cite{Otal})}\\
\label{thm: Otal}
There exists a positive constant $c_0$ such that
for a marked quasifuchsian punctured torus group $G(q)$
and $V \in G(q)$ representing a simple closed geodesic $\gamma$  in $\HH^3/G(q)$,
if  $V$ is purely hyperbolic and the hyperbolic length $\l_{\gamma}(q)$ of
$\gamma$ is less than $c_0$,
then $\gamma$ is a bending locus of $\bch(G(q))$.
\end{thm}

Following the proof of proposition 9 in~\cite{Otal}, we will give a proof of
theorem~\ref{thm: Otal}
to estimate $c_0$ in the next subsection 5.2.

Suppose that
$\gamma$ is not the bending locus of $\bch(G(q))$.
Then the pleated surface $\psi$ is not a homeomorphism by
proposition~\ref{prop: homeo}.
Let $H$ be a Fuchsian subgroup $\langle V, WVW^{-1} \rangle$
where $(V, W)$ is a marking of $G(q)$.
Denote the totally geodesic plane whose boundary $\partial \DD$ contains
$\Lambda(H)$ by $\DD \subset \HH^3$.
Let $P$ be the convex hull of $\Lambda(H)$ in $\DD \subset \HH^3$.
Then $g(P)$ is the convex hull of $\Lambda(gHg^{-1})$ in $g(\DD) \subset
\HH^3$.
Now we have a following claim.

\begin{prop}
\label{prop: intersection}
If the pleated surface $\psi$ is not a homeomorphism,
then there exists $g \in G(q)$ such that $P$ and $g(P)$ intersect
transversally in the axis of $V$.
\end{prop}

To show this proposition, we need two lemmas.

\begin{lem}
\label{lem: group}
For $g \in G(q)$, $gHg^{-1} \cap H$ is trivial or cyclic subgroup generated
by $gVg^{-1}$.
\qed
\end{lem}

\begin{lem}
\label{lem: limit set}
For $g \in G(q)$,
$\Lambda(gHg^{-1}) \cap \Lambda(H)$ is empty or fixed points of $gVg^{-1}$.
\end{lem}
\proof
$H$ and $gHg^{-1}$ are Fuchsian subgroups of a quasifuchsian group $G(q)$.
From a theorem of Suskind (see theorem 3.14 in~\cite{MT}),
$$
\Lambda(gHg^{-1}) \cap \Lambda(H)= \Lambda(gHg^{-1} \cap H).
$$
Hence it concludes the proof by using lemma~\ref{lem: group}.
\qed

\vspace{3mm}

Now we can show proposition~\ref{prop: intersection}.
Since we assume that $\psi$ is not a homeomorphism, there exists $g \in G$
such that the interior of $P$ and the interior of $g(P)$ intersect
transversally.
Then from lemma~\ref{lem: limit set}, $\Lambda(gHg^{-1}) \cap \Lambda(H)$ is
empty.
Therefore the axis of $V$ cuts $g(P)$ transversally.
\qed

\vspace{3mm}

Now we assume that $c_0$ is smaller than the Margulis constant.
Then the interior of $g(P)$ cut the Margulis tube $T$ with radius $r$ along
the axis of $V$ transversally. Hence now we have a geodesic disk $\Delta
= T \cap g(P)$ on $g(P)$.

\begin{lem}
The hyperbolic area of $\Delta$ is bigger than $4 \pi \sinh^2 (r/2)$.
\end{lem}
\proof
If $g(P)$ intersects the axis of $V$ orthogonally, then $\Delta$ is a
hyperbolic disk of radius $r$, hence
the hyperbolic area of it is $4 \pi \sinh^2 (r/2)$.
If $g(P)$ intersects the axis of $V$ not orthogonally, then $\Delta$
contains a hyperbolic disk of radius $r$, hence
the hyperbolic area of $\Delta$ is bigger than $4 \pi \sinh^2 (r/2)$.
\qed

\vspace{3mm}

Now we can give a proof of theorem~\ref{thm: Otal}.
By the Margulis lemma, $\Delta$ projects into the image of $g(P)$ in
$\HH^3/G(q)$ injectively,
whereas the image of $g(P)$ in  $\HH^3/G$ has its hyperbolic area  $2 \pi$
since it is the isometric
image of a punctured cylinder.
Therefore if the hyperbolic length of $\gamma$ in $\HH^3/G(q)$ is
sufficiently small, we can take a radius $r$ of the Margulis tube $T$ satisfying $4 \pi
\sinh^2 (r/2) \geq 2 \pi$, which is a contradiction.
This concludes the theorem.
\qed

\subsection{A lower bound of $c_0$}
\label{sec: lower}
Following~\cite{Meyer}, we have a formula of the radius of a Margulis tube.

\begin{prop}{\rm (c.f. theorem in section 3 of~\cite{Meyer})}\\
For $q \in \QF$,
assume that $V \in G(q)$ representing a simple closed geodesic $\gamma$ in $\HH^3/G(q)$, 
which is purely hyperbolic.
If the hyperbolic length $\l_{\gamma}$ of $\gamma$ satisfies $\cosh \l_{\gamma}< \sqrt{2}$,
then there is a Margulis tube with radius $r$ satisfying
$$
\sinh^2(r)=\frac{1}{2}(\frac{\sqrt{3-2\cosh \l_{\gamma}}}{\cosh
\l_{\gamma}-1}-1).
$$
\qed
\end{prop}
The inequality $4 \pi \sinh^2 (r/2)  \geq 2 \pi$ and the above formula
give us a lower bound of $c_0$.
\begin{cor}
$\cosh^{-1} \frac{48+5\sqrt{2}}{49} \approx 0.493 \leq c_0$.
\qed
\end{cor}

\section{Existence of non-standard components}

\begin{thm}
\label{thm: existence}
There is some positive constant $c_1$ such that for any $c$ satisfying $c>
c_1$,
$\LL$ contains non-standard components.
\end{thm}
To prove this theorem, we use the Earle slice $\E$
of punctured tori studied in~\cite{Kom, KoS}.
This idea is due to Raquel Diaz.
We review notations of $\E$ (c.f.~\cite{KoS}).
The {\it Earle slice} $\E$ of $\QF$ is the set of
$G(q)=\langle A,B \rangle$ satisfying the following symmetry;
there exists an elliptic element of
order 2 such that $EAE=B$.
Then $\E$ is a holomorphic slice of $\QF$ and considering the
conformal structure of $\Omega_+(G(q))/G(q)$, it is naturally isomorphic to
the Teichm\"{u}ller
space of punctured tori.
Any element of $\E$ can be represented by the following matrices
in $SL(2,\CC)$ of the form
$A= A_d, B= B_d , d \in \CC - \{0\}$,  where

$$
A_d  =  \left(
\begin{array}{cc}
\frac{d^2+1}{d} & \frac{d^3}{2d^2+1}\\
\frac{2d^2+1}{d} & d
\end{array}
\right) ,
B_d  =  \left(
\begin{array}{cc}
\frac{d^2+1}{d} & -\frac{d^3}{2d^2+1}\\
-\frac{2d^2+1}{d} & d
\end{array}
\right).
$$
The complex parameter $d$ gives a holomorphic embedding of $\E$
into the right half plane $\CC^+$ and we assume that $\E$ is embedded in
$\CC^+$.
Then $\E$ contains the positive real line $\RR^+$ which is the Fuchsian
locus $\E \cap \F$ of $\E$.
Put ${\bf C}^+_{d'}=\{ d \in {\bf C}^+ | \Re \; d > d' \}$.
To show theorem~\ref{thm: existence}, we need lemmas.

\begin{lem}
There is a positive constant $d_0$ such that for any $d'>d_0$,
the hyperbolic locus $\H_{2/1}$ of $T_{2/1}(d) = \tr W_{2/1}(d)$
satisfies
$$
(\H_{2/1}- {\bf R}^+) \cap {\bf C}^+_{d'}
\neq \emptyset.
$$
\end{lem}
\proof
We remark that $W_{2/1}=A^{-2}B$.
Then we can check our claim by direct calculation.
\qed

\begin{lem}
There is a positive constant $d_1$ such that 
the $2/1$-pleating variety $\P_{2/1}$ satisfies
$$
\P_{2/1} \cap \E \subset {\bf C}^+ - {\bf C}^+_{d_1}.
$$
\end{lem}
\proof
In~\cite{KoS}, it is shown that $\P_{p/q} \cap \E$ is equal to two
components of $\H_{p/q}- {\bf R}^+$
terminating to the unique critical point of $\tr W_{p/q}$ on ${\bf R}^+$
(c.f. theorem 5.1 in~\cite{KoS}).
Then we can check our claim by direct calculation.
\qed

\begin{lem}{\rm (c.f.~\cite{Kom})}
There is a positive constant $d_2$ such that
$$
{\bf C}^+_{d_2} \subset \E.
$$
\qed
\end{lem}

Now we can prove  theorem~\ref{thm: existence}.
There is a positive constant $c_1$ such that for any
$c > c_1$, there is $d \in \E$ such that
the word $A^{-2}B$ is purely hyperbolic and $\l_{W_{2/1}}(d)=c$,
but $d$ is not contained in $\P_{2/1}$.
This concludes the theorem.
\qed

\begin{remark}
To estimate $c_1$, we need to know the size of the round disk contained in $\E$
tangent to the boundary $\partial \E$ of $\E$ at the origin (see~\cite{Kom}).
\end{remark}

\vspace{3mm}

Comparing with the results in section 5 and 6, 
we have the following conjecture supported by numerical experiences 
by the second author.

\begin{conj}
There exists a unique $c_0$ such that
$\LL$ coincides with the standard component for any $c \leq c_0$,
while
$\LL$ contains infinitely many non-standard components for any $c>c_0$.
\end{conj}

\section{Scaling property of $\LL$}

In the final section we will study the self-similar phenomena of $\LL$ which we can observe from 
figures of  $\LL$ in this paper. 
First we remark that  $\LL$ has analytic automorphisms coming from Dehn twists.

\begin{prop}
$(A,B) \in \LL$ implies $(A,A^nB) \in \LL$ for all $n \in \ZZ$.
\end{prop}

\proof
The automorphism of $G$ defined by $(A,B) \mapsto (A,A^nB)$ is a
Dehn twist along $A$ which preserves $\QF$ and $\Q_{\gamma, c}$.
\qed

\vspace{3mm}

Next result is easy to prove, but it induces the asymptotic self similarity
of $\LL$.
\begin{prop}
$$
\lim_{n \rightarrow \infty}
\frac{\tr A^nB}{\tr A^{n-1}B}
= \tr A \cdot
\frac{1+ \sqrt{1-(\frac{2}{\tr A})^2}}{2},
$$
which is the attractive fixed point of the map
$\tr A - \frac{1}{z}$.
\end{prop}

\proof
The following trace identity is well known:
$$
\tr A^nB = \tr A \cdot \tr A^{n-1}B - \tr A^{n-2}B.
$$
Divide both sides by $(\tr A)^n$ and put
$x_n:=\frac{\tr A^nB}{(\tr A)^n}$.
Then we have
$$
x_n=x_{n-1}-\frac{1}{(\tr A)^2}x_{n-2}.
$$
Moreover put
$y_n:=\frac{x_n}{x_{n-1}}$, then
$$
y_n=1-\frac{1}{(\tr A)^2}\frac{1}{y_{n-1}}.
$$
Finally put
$z_n:=\tr A \cdot y_n$, then
$$
z_n=\tr A-\frac{1}{z_{n-1}}.
$$
Since  $A$ is purely hyperbolic, the linear fractional transformation
$$
w=Tr A-\frac{1}{z}
$$
is also purely hyperbolic, hence
all points besides the repelling fixed point of $A$
converge to the attracting fixed point of $A$,
$\tr A \cdot \frac{1+ \sqrt{1-(\frac{2}{\tr A})^2}}{2}$.
From the above arguments, $z_n=\frac{\tr A^nB}{\tr A^{n-1}B}$ converges to
this point.
\qed

\begin{cor}
Linear slice has an asymptotic scaling constant
$\tr A \cdot \frac{1+ \sqrt{1-(\frac{2}{\tr A})^2}}{2}$.
\qed
\end{cor}

\begin{rem}
When $A$ tends to be parabolic,
$$
\lim_{n \rightarrow \infty}
\frac{\tr A^nB}{\tr A^{n-1}B}
=1,
$$
which relates to the fact that the Maskit slice is invariant
under translations.
\end{rem}

\begin{figure}
\begin{center}
\includegraphics[width=5cm]{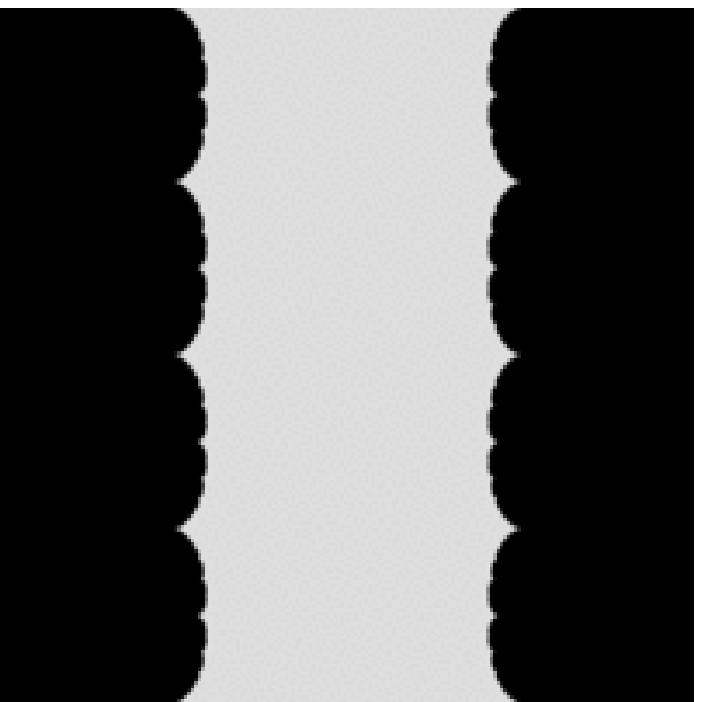}
\includegraphics[width=5cm]{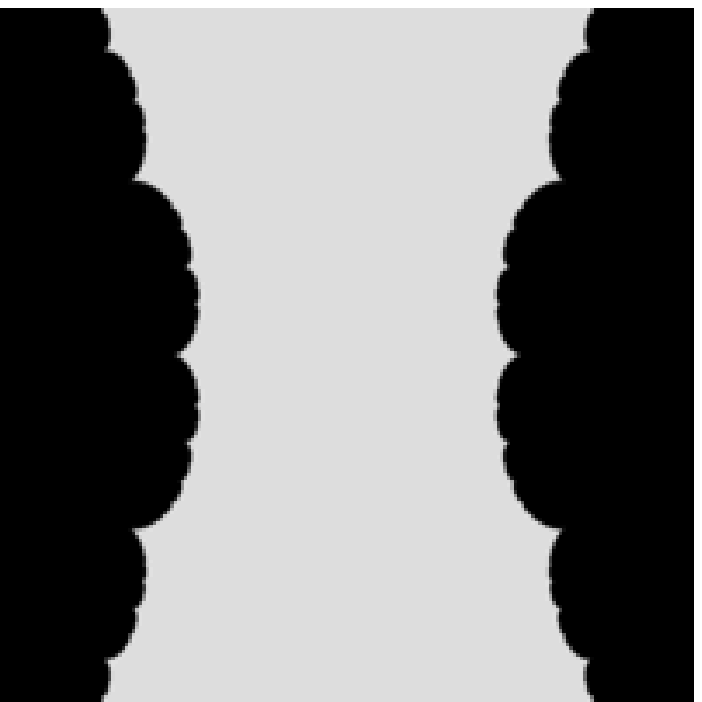}
\caption{Maskit slice (left) and $\tr A=2.5$ slice (right)}
\end{center}
\end{figure}

\begin{figure}
\begin{center}
\includegraphics[width=3.5cm]{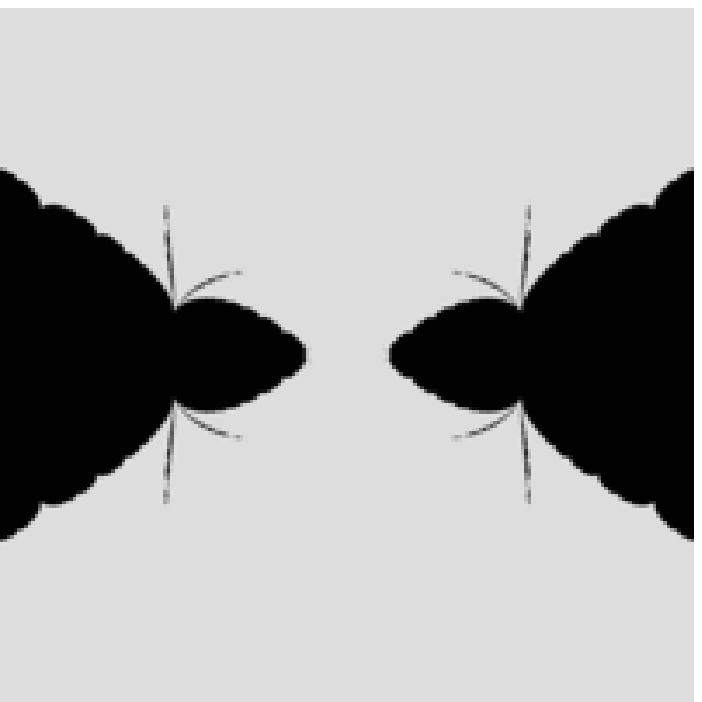} \
\includegraphics[width=3.5cm]{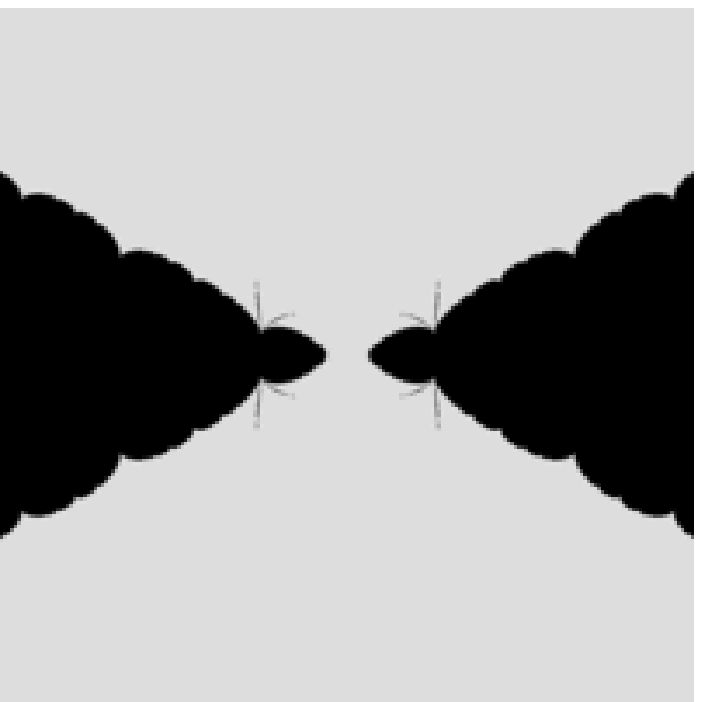} \
\includegraphics[width=3.5cm]{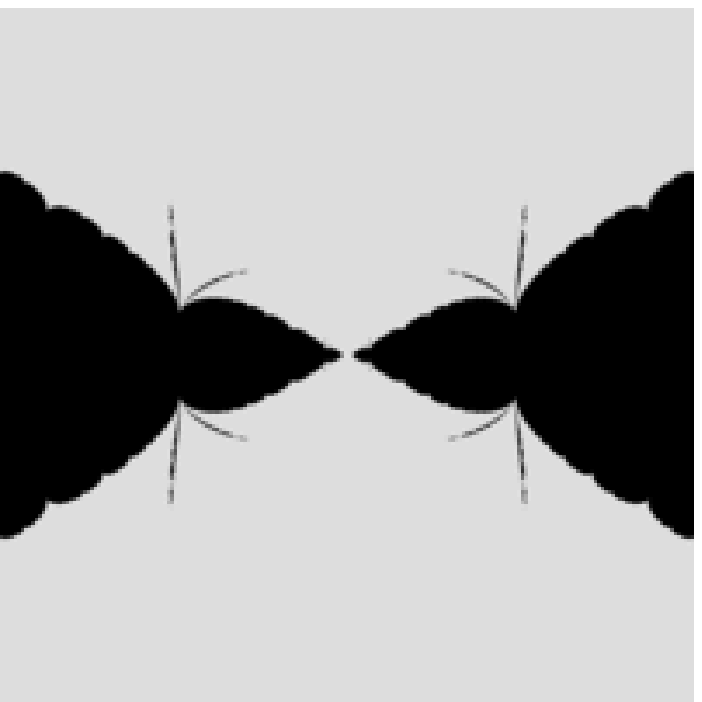}
\caption{$\tr A=8$ linear slices with ranges 16(left), 32(center), 128(right)}
\end{center}
\end{figure}

\begin{figure}
\begin{center}
\includegraphics[width=3.5cm]{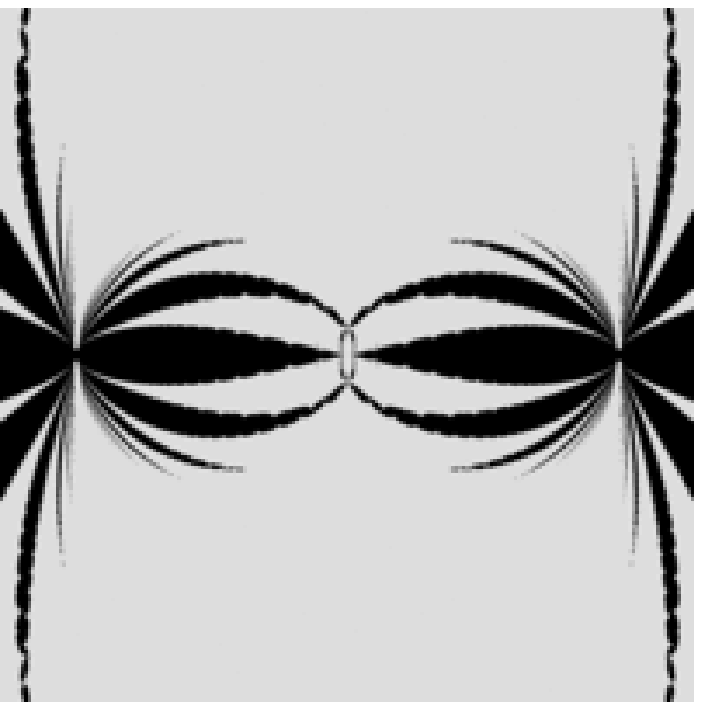} \
\includegraphics[width=3.5cm]{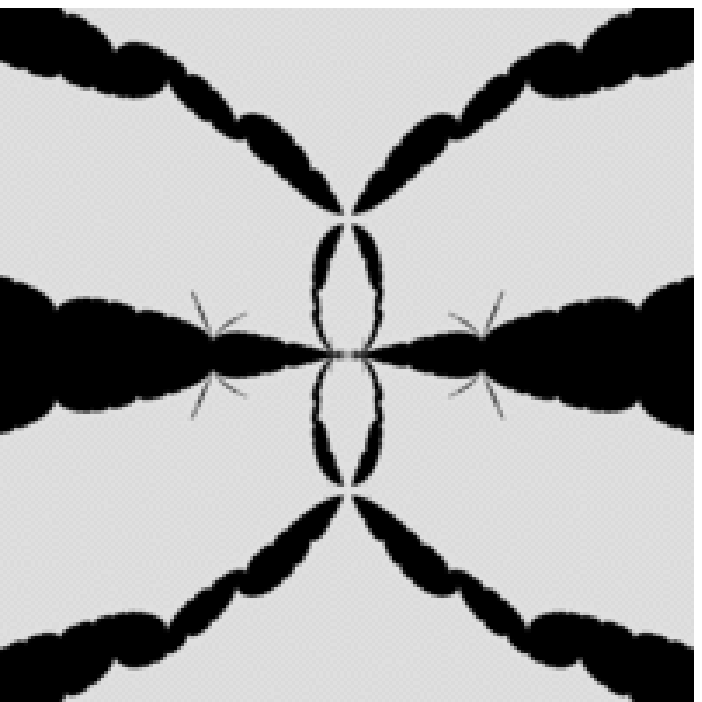} \
\includegraphics[width=3.5cm]{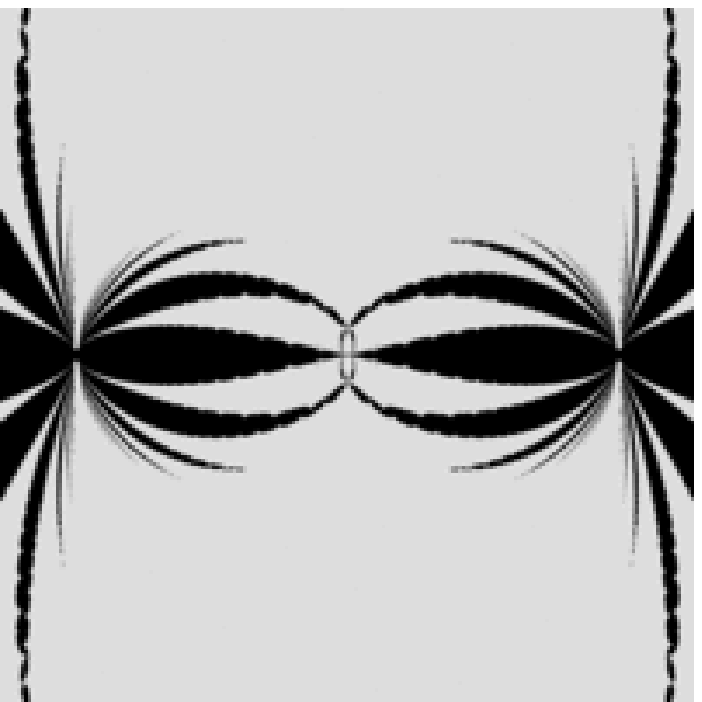}
\caption{$\tr A=100$ linear slices with ranges 128(left), 2560(center), 12800(right)}
\end{center}
\end{figure}

\begin{rem}
Even if $A$ is loxodromic, $\LL$ has this scaling property.
Hence we can also see that the figure 10 in~\cite{Mc}
also has such scaling property.
\end{rem}

\end{document}